\documentstyle[11pt,amstex,amssymb]{amsart}\textheight=19.8truecm\textwidth=13.7truecm
\begin{document} 

\vskip 30pt
\centerline{\Large\bf A concavity inequality for symmetric norms}
\vskip 30pt
\centerline{Jean-Christophe Bourin}
\vskip 10pt
 \centerline{E-mail: bourinjc@@club-internet.fr}
 \vskip 5pt
 \centerline{Universit\'e de Cergy-Pontoise, d\'ept.\ de Math\'ematiques}
 \vskip 5pt
 \centerline{2 rue Adolphe Chauvin, 95302 Pontoise, France}

\vskip 20pt
\noindent
{\small {\bf Abstract.} 
\vskip 5pt
We review  some recent convexity results for Hermitian matrices and we add a new one to the list: Let  $A$ be semidefinite positive, let $Z$ be expansive, $Z^*Z\ge I$, and  let $f:[0,\infty)\longrightarrow[0,\infty)$ be a concave function. Then, for all symmetric norms
$$
\Vert f(Z^*AZ)\Vert \le \Vert Z^*f(A)Z\Vert.
$$
This inequality complements a classical trace inequality of Brown-Kosaki.


Mathematical subjects classification:   47A30 47A63}

\vskip 25pt
{\large\bf Introduction}
\vskip 10pt
A good part of Matrix Analysis consists in establishing results for Hermitian operators considered as generalized real numbers. In particular several results are
 matrix versions of inequalities for convex functions $f$ on the real line, such as
\begin{equation}
f\left(\frac{a+b}{2}\right) \le \frac{f(a)+f(b)}2
\end{equation}
for all reals $a$, $b$ and
\begin{equation}
f(za)\le zf(a)
\end{equation}
for convex functions $f$ with $f(0)\le 0$ and scalars $a$ and  $z$ with $0<z<1$.

In this brief note we first review some recent matrix versions of (1), (2) and next we give the matrix version of the companion inequality of (2):
\begin{equation}
f(za)\le zf(a)
\end{equation}
for concave functions $f$ with $f(0)\ge 0$ and scalars $a$ and  $z$ with $1<z$.

Capital letters $A$, $B\dots,Z$ mean $n$-by-$n$ complex matrices, or operators on a finite dimensional Hilbert space ${\cal H}$; $I$ stands for the identity. When $A$ is positive semidefinite, resp.\ positive definite, we write $A\ge 0$, resp.\ $A>0$. 

\vskip 40pt\noindent
{\large\bf 1. Some known convexity results }
\vskip 10pt
The following are wellknown trace versions of elementary inequalities (1) and (2).

\vskip 10pt\noindent 1.1. von Neuman's Trace Inequality: For convex functions $f$ and Hermitians $A$, $B$,
\begin{equation}
{\rm Tr}\,f\left(\frac{A+B}{2}\right) \le {\rm Tr}\,\frac{f(A)+f(B)}{2}
\end{equation}
equivalently ${\rm Tr}\circ f$ is convex on the set of Hermitians.

\vskip 10pt\noindent
1.2. Brown-Kosaki's Trace Inequality [5]: Let $f$ be convex with  $f(0)\le 0$ and let $A$ be Hermitian. Then, for  all contractions $Z$, 
\begin{equation}
{\rm Tr}\,f(Z^*AZ) \le {\rm Tr}\,Z^*f(A)Z. 
\end{equation}

\vskip 10pt\noindent
1.3.  Hansen-Pedersen's Trace Inequality [7]:    Let $f$ be convex  and let $\{A_i\}_{i=1}^n$ be Hermitians.  Then, for  all isometric columns $\{Z_i\}_{i=1}^n$,
\begin{equation*}
{\rm Tr}\,f(\sum_iZ_i^*A_iZ_i) \le {\rm Tr}\,\sum_iZ_i^*f(A_i)Z_i. 
\end{equation*}

\vskip 10pt\noindent
Here isometric column means that $\sum_i Z_i^*Z_i=I$. Hansen-Pedersen's result contains (4) and (5).

\vskip 10pt 
When $f$ is convex and monotone, we showed [2] that the above trace inequalities can be extended to  operator inequalities up to a unitary congruence. Equivalently we have inequalities for eigenvalues. Let us give the precise statements corresponding to von Neumann and Brown-Kosaki trace inequalities.

\vskip 10pt\noindent
1.4. Let $A$, $B$ be Hermitians and let $f$ be a $monotone$ convex function. Then,  there exists a unitary $U$ such that
\begin{equation}
f\left(\frac{A+B}{2}\right) \le U\cdot\frac{f(A)+f(B)}{2}\cdot U^*
\end{equation}

\vskip 10pt\noindent
1.5. Let $A$ be a  Hermitian, let $Z$ be  a contraction and let $f$ be a $monotone$ convex function. Then, there exists a unitary $U$ such that
\begin{equation}
f(Z^*AZ) \le UZ^*f(A)Z U^*
\end{equation}

\vskip 10pt
Statements 1.4 and 1.5 can break down when the monotony assumption is dropped. But we recently obtained [4] substitutes
involving the mean of two unitary congruences. Let us recall the precise result corresponding to inequalities (1) and (6). 

\vskip 10pt\noindent
1.6. Let $f$ be a convex function,  let $A$, $B$ be Hermitians and set $X=f(\{A+B\}/2)$ and $Y= \{f(A)+f(B)\}/2$. Then, there exist  unitaries  $U$, $V$ such that
\begin{equation*}
X \le \frac{UYU^* + VYV^*}{2}.
\end{equation*}

\vskip 10pt
Another substitute of (6) for general convex functions $f$ would be a positive answer to the following still open problem [2]: 
Given Hermitians $A$, $B$, can we find unitaries $U$, $V$ such that
\begin{equation*}
f\left(\frac{A+B}{2}\right) \le \,\frac{Uf(A)U^*+Vf(B)V^*}{2} \ ?
\end{equation*}

\vskip 10pt
We turn to a Brown Kosaki type inequality involving expansive operators $Z$, that is $Z^*Z\ge I$. We showed the following trace version of the elementary inequality (3):

\vskip 10pt\noindent
1.7.  Let $f$ be convex with  $f(0)\le 0$ and let $A\ge 0$. Then, for  all expansive operators $Z$, 
\begin{equation}
{\rm Tr}\,f(Z^*AZ) \ge {\rm Tr}\,Z^*f(A)Z. 
\end{equation}

\vskip 10pt
It is interesting to note [2] that, contrarily to the contractive case (5), the assumption $A\ge 0$ can not be dropped.
Also, still contrarily to (5), this result can not be extended to eigenvalues inequalities like (7). Nevertheless, we have:

\vskip 10pt\noindent
1.8.  Let $f$ be nonnegative convex with  $f(0)=0$, let $A\ge 0$ and let $Z$ be  expansive. Then, for all symmetric norms 
\begin{equation}
\Vert f(Z^*AZ)\Vert \ge \Vert Z^*f(A)Z\Vert. 
\end{equation}

\vskip 10pt\noindent
Here, by symmetric norm we mean a unitarily invariant one, that is $\Vert A\Vert= \Vert UAV\Vert$ for all operators $A$ and all unitaries $U$, $V$.

\vskip 20 pt\noindent
{\large\bf 2. A new concavity result}
\vskip 10pt
Of course if $f$ is concave with $f(0)\ge 0$ then inequality (8) is reversed and provides an extension of its scalar version (3). Assuming furthermore $f$ nonnegative we tried to extend it to all symmetric norms but, besides the trace norm, we  only got the operator norm case. Here we may state: 

\vskip 10pt \noindent
{\bf Theorem 2.1.} {\it Let $f:[0,\infty)\longrightarrow [0,\infty)$ be a concave function. Let $A\ge0$  and let $Z$ be expansive. Then, for all symmetric norms
\begin{equation*}
\Vert f(Z^*AZ) \Vert \le \Vert Z^*f(A)Z\Vert. 
\end{equation*}
}

\vskip 10pt\noindent
{\bf Proof.} It suffices to prove the theorem for the Ky Fan $k$-norms $\Vert\cdot\Vert_{k}$ (cf.\ [1]). This shows, since $Z$ is expansive, that we may assume that $f(0)=0$. Note that $f$ is necessarily nondecreasing. Hence, there exists a rank $k$ spectral projection $E$ for $Z^*AZ$, corresponding to the $k$-largest eigenvalues $\lambda_1(Z^*AZ),\dots,\lambda_k(Z^*AZ)$ of $Z^*AZ$, such that
$$
\Vert f(Z^*AZ)\Vert_{k}=\sum_{j=1}^k \lambda_j(Z^*AZ)={\rm Tr\,} Ef(Z^*AZ)E.
$$
Therefore, using a wellknown property of Ky Fan norms, it suffices to show that
$$
{\rm Tr\,} Ef(Z^*AZ)E \le {\rm Tr\,} EZ^*f(A)ZE.
$$
This is the same as requiring that
\begin{equation}
{\rm Tr\,} EZ^*g(A)ZE \le {\rm Tr\,} Eg(Z^*AZ)E
\end{equation}
for all convex functions $g$ on $[0,\infty)$ with $g(0)=0$. Any such function can be approached by a combination of the type
\begin{equation}
g(t)=\lambda t +\sum_{i=1}^n \alpha_i(t-\beta_i)_+
\end{equation}
for a scalar $\lambda$ and some nonnegative scalars $\alpha_i$ and $\beta_i$. Here $(x)_+=\max\{0,x\}$. By using the linearity of the trace it suffices to show that (10) holds for $g_{\beta}(t)=(t-\beta)_+$, $\beta\ge 0$. We claim that there exists a unitary $U$ such that
\begin{equation}
Z^*g_{\beta}(A)Z \le Ug_{\beta}(Z^*AZ)U^*.
\end{equation}
This claim and a basic property of the trace then show that (10) holds for $g_{\beta}$. Indeed, we then have
\begin{align*}
{\rm Tr\,} EZ^*g_{\beta}(A)ZE &= \sum_{j=1}^k \lambda_j(EZ^*g_{\beta}(A)ZE) \\
&\le \sum_{j=1}^k \lambda_j(Z^*g_{\beta}(A)Z) \\
&\le \sum_{j=1}^k \lambda_j(g_{\beta}(Z^*AZ)) \quad {\rm (by\ 12)} \\
&= \sum_{j=1}^k \lambda_j(Eg_{\beta}(Z^*AZ)E) \\
&= {\rm Tr\,} Eg_{\beta}(Z^*AZ)E 
\end{align*}
where the fourth equality follows from the fact that $g_{\beta}$ is nondecreasing
and hence $E$ is also a spectral projection of $g_{\beta}(Z^*AZ)$ corresponding to the $k$ largest eigenvalues.

The inequality (12) has been established in [2] in order to prove (8). Let us recall the proof of (12):
We will use the following simple fact. If $B$ is a positive operator with ${\rm Sp}B\subset\{0\}\cup(x,\infty)$, then we also have ${\rm Sp}Z^*BZ\subset\{0\}\cup(x,\infty)$. Indeed $Z^*BZ$ and $B^{1/2}ZZ^*B^{1/2}$ (which is greater than $B$) have the same spectrum.

Let $P$ be the spectral projection of $A$ corresponding to the eigenvalues strictly greater than $\beta$
and let $A_{\beta}=AP$. Since $Z^*AZ-\beta I\ge Z^*A_{\beta}Z-\beta I$ and $t\longrightarrow t_+$ is nondecreasing, there exists a unitary operator $V$ such that 
\begin{equation*}
 (Z^*AZ-\beta I)_+\ge V(Z^*A_{\beta}Z-\beta I)_+V^*
\end{equation*}
Since $Z^*(A-\beta I)_+Z=Z^*(A_{\beta}-\beta I)_+Z$ we may then assume that $A=A_{\beta}$. Now, the above simple fact implies
\begin{equation*}
 (Z^*A_{\beta}Z-\beta I)_+=Z^*A_{\beta}Z-\beta Q
\end{equation*}
where $Q={\rm supp}Z^*A_{\beta}Z$ is the support projection of $Z^*A_{\beta}Z$. Therefore, using $(A_{\beta}-\beta I)_+=A_{\beta}-\beta P$, it suffices to show the existence of a unitary operator $W$ such that
\begin{equation*}
 Z^*A_{\beta}Z-\beta Q  \ge WZ^*(A_{\beta}-\beta P)ZW^*=WZ^*A_{\beta}ZW^*-\beta WZ^*PZW^*.
\end{equation*}
But, here we can take $W=I$. Indeed, we have
\begin{equation*}
 {\rm supp}Z^*PZ=Q \ (*) \quad{\rm and}\quad {\rm Sp}Z^*PZ\subset\{0\}\cup[1,\infty) \ (**)
\end{equation*}
where ($\ast\ast$) follows from the above simple fact and the identity ($\ast$) from the  observation below with $X=P$ and $Y=A_{\beta}$.
\vskip 5pt\noindent
{\it Observation.}\, If $X$, $Y$ are two positive operators with ${\rm supp}X={\rm supp}Y$, then for every  operator $Z$ we also have ${\rm supp}Z^*XZ={\rm supp}Z^*YZ$. 
\vskip 5pt\noindent
To check this, we establish the corresponding equality for the kernels,
\begin{equation*}
\ \qquad \ker Z^*XZ=\{h\ :\ Zh\in \ker X^{1/2}\} = \{h\ :\ Zh\in \ker Y^{1/2}\}= \ker  Z^*YZ. \qquad \Box
\end{equation*}

\vskip 10pt
In the above proof, the simple idea of approaching convex functions as in (11) was fruitful. It is also useful to prove (see [2]) the Rotfel'd Trace Inequality: For  concave functions $f$ with $f(0)\ge0$  and $A$, $B\ge 0$,
$$
{\rm Tr\,}f(A+B)\le {\rm Tr\,}f(A)+{\rm Tr\,}f(B).
$$
If $f$ is convex with $f(0)\le0$ the reverse inequality holds, in particular we have
McCarthy's inequality
$$
{\rm Tr\,}(A+B)^p\ge {\rm Tr\,}A^p+{\rm Tr\,}B^p
$$
for all $p>1$.

\vskip 10pt\noindent
{\bf Remark 2.2.} Though scalars inequalities (2), (3) or their concave anologous hold for a more general class than convex or concave functions,
the corresponding trace inequalities need the convexity or concavity assumption (cf.\ [2]). A fortiori, Theorem 2.1 needs the concavity assumption.

\vskip 10pt\noindent
{\bf Remark 2.3.} When $f$ is operator monotone, Theorem 2.1 extends to an operator inequality which can be rephrased for contractions  as follows: For nonnegative operator monotone functions $f$ on $[0,\infty)$, contractions $Z$  and $A\ge 0$,
$$
Z^*f(A)Z \le f(Z^*AZ)
$$
This is the famous Hansen's inequality [4]. Similarly when $f$ is operator convex, Hansen-Pedersen's Trace Inequality can be extended to an operator inequality [7] (see also [3]).

\vskip 10pt
 Extensions of Theorem 2.1  to infinite dimensional spaces will be considered in a forthcoming work.
 
\vskip 10pt\noindent
{\bf Acknowledgement.} The author is grateful to a referee for its valuable comments.

\vskip 10pt
{\bf References}

\noindent
{\small 

\noindent
[1] R.\ Bhatia, Matrix Analysis, Springer, Germany, 1996.

\noindent
[2] J.-C.\ Bourin, {\it Convexity or concavity inequalities for Hermitian operators,}   Math.\ Ineq.\ Appl. 7 ${\rm n}^04$  (2004) 607-620.

\noindent
[3] J.-C.\ Bourin, Compressions, Dilations and Matrix Inequalities, RGMIA monograph, Victoria university, 
Melbourne 2004 (http://rgmia.vu.edu.au/monograph)

\noindent
[4] J.-C.\ Bourin, {\it  Hermitian operators and convex functions}, to appear in J.\  Math.\ Ineq.

\noindent
[5] L.\ G.\ Brown and H.\ Kosaki, {\it Jensen's inequality in semi-finite von Neuman algebras}, J. Operator Theory 23 (1990) 3-19.

\noindent
[6] F.\ Hansen, {\it An operator inequality}, Math.\ Ann.\ 246 (1980) 249-250.

\noindent
[7] F.\ Hansen and G.\ K.\ Pedersen, {\it Jensen's operator inequality}, Bull. London Math.\ Soc.\
35 (2003) 553-564.

\end{document}